\documentclass[12pt]{article}
\usepackage[latin1]{inputenc}
\usepackage{amsmath}
\usepackage{amsfonts}
\usepackage{amssymb}
\usepackage{amsthm}

\usepackage{tikz}

\usepackage{fixltx2e} % subscripts

\usepackage{titlesec}
\titleformat{\subsection}[hang]{\normalfont\bfseries}{\thesubsection}{1em}{}
\titlespacing\section{0pt}{3.5ex plus 0.5ex minus .2ex}{0.3ex plus .2ex}
\titlespacing\subsection{0pt}{2.5ex plus 0.5ex minus .2ex}{0.3ex plus .2ex}
\titlespacing\subsubsection{0pt}{2.5ex plus 0.5ex minus .2ex}{0.3ex plus .2ex}
\usepackage{mathrsfs} % for \mathscr (script letters)

\usepackage[alphabetic,lite]{amsrefs} % for different reference style

\usepackage{fullpage}
\usepackage{setspace}
\usepackage{hyperref}
\usepackage{color}
\usepackage{enumitem}
\usepackage{ulem} %um text durchzustreichen

\usepackage[all,cmtip]{xy} % crazy diagrams

 \usepackage{fancyhdr}
\newtheorem{Thm}{Theorem}[section]

\newtheorem{Lemma}[Thm]{Lemma}

\newtheorem{Cor}[Thm]{Corollary}

\newtheorem{Prop}[Thm]{Proposition}

\theoremstyle{definition}
\newtheorem{Def}[Thm]{Definition}
\newtheorem{Rem}[Thm]{Remark}

\newcommand{\Proof}{\textbf{Proof. }}

\newcommand{\<}{\left\langle}
\renewcommand{\>}{\right\rangle}
\newcommand{\bG}{\mathbb{G}}

\newcommand{\bN}{\mathbb{N}}

\newcommand{\bQ}{\mathbb{Q}}
\newcommand{\bR}{\mathbb{R}}

\newcommand{\bZ}{\mathbb{Z}}

\newcommand{\cO}{\mathcal{O}}
\newcommand{\cP}{\mathcal{P}}

\newcommand{\cS}{\mathcal{S}}

\newcommand{\cW}{\mathcal{W}}

\newcommand{\fg}{\mathfrak{g}}

\newcommand{\ft}{\mathfrak{t}}

\DeclareMathAlphabet{\mathpzc}{OT1}{pzc}{m}{it}

\newcommand{\ra}{\rightarrow}

\newcommand{\ov}{\overline}

\newcommand{\ul}{\underline}

\providecommand{\abs}[1]{\left\lvert#1\right\rvert}

\DeclareMathOperator{\Hom}{Hom}			% Set of arrows between two object
\DeclareMathOperator{\SL}{SL}		  	% SL
\DeclareMathOperator{\Dyn}{Dyn}			% Dynkin diagram
\DeclareMathOperator{\Gal}{Gal}			% Galois group
\DeclareMathOperator{\Aut}{Aut}			% automorhpism group
\DeclareMathOperator{\Lie}{Lie}  	  % Lie algebra
\DeclareMathOperator{\Res}{Res}  	  % restriction
\DeclareMathOperator{\val}{val}  	  % valuation
\DeclareMathOperator{\Cent}{Cent}  	  % centralizer
\newcommand{\ba}{\begin{aligned}}
\newcommand{\ea}{\end{aligned}}
\newcommand{\beqn}{\begin{eqnarray}}
\newcommand{\eeqn}{\end{eqnarray}}
\newcommand{\beqns}{\begin{eqnarray*}}
\newcommand{\eeqns}{\end{eqnarray*}}
\newcommand{\benum}{\begin{enumerate}}
\newcommand{\eenum}{\end{enumerate}}

\newcommand{\Wabs}{\cW}

\def\Jtype/{$J$-typ\ul{e}}
\def\almoststable/{almost stabl\ul{e}}
\def\pure/{pur\ul{e}}
\def\dualdepth/{dual dept\ul{h}}
\def\nonamesequence/{noname sequenc\ul{e}}
\def\extendeddatum/{extended datu\ul{m}}
\def\datum/{datu\ul{m}}
\def\truncatedextendeddatum/{truncate\ul{d} extended datu\ul{m}}
\def\depthpreserving/{depth preservin\ul{g}}
\def\wellbehaved/{well behave\ul{d}}
\def\generic/{generi\ul{c}}
\def\tame/{tam\ul{e}}

\AtEndDocument{\bigskip{\footnotesize%
		\textsc{IAS Fuld Hall 408, 1 Einstein Drive, Princeton, NJ 08540, USA} \par  
		\textit{E-mail address}: \texttt{fintzen@umich.edu}
	}}

\begin{document}
\author{J. Fintzen\\[1cm]}
\title{Tame tori in $p$-adic groups and good semisimple elements}
\date{}
\maketitle

	\begin{abstract} 
	Let $G$ be a reductive group over a non-archimedean local field $k$. We provide necessary conditions and sufficient conditions for all tori of $G$ to split over a tamely ramified extension of $k$. We then show the existence of good semisimple elements in every Moy--Prasad filtration coset of the group $G(k)$ and its Lie algebra, assuming the above sufficient conditions are met. \\[1cm]
	
	\end{abstract}

	{
		\renewcommand{\thefootnote}{}  % to delete the footnote number 	
		\footnotetext{MSC2010: Primary 20G25, 20G07; Secondary 22E50, 22E35, 20C08} % 22E50 (reps of lin alg groups over local fields), 20G25 (red groups over local fields), 11S37 (LLC, I would say secondary?), 14L24 (GIT, maybe secondary?); 20G15: 	Linear algebraic groups over arbitrary fields;  	
		%20C08 Hecke algebras 
		%20G07 structure theory  20G05 rep thy
		%  	22E35   	Analysis on $p$-adic Lie groups
		\footnotetext{Keywords: reductive groups over non-archimedean local fields, tamely ramified maximal tori, order of the Weyl group, Moy--Prasad filtrations, good semisimple elements}
	}

\tableofcontents

\newpage
\section{Introduction}
Many important results in the representation theory of  %and harmonic analysis on
 reductive groups over non-archimedean local fields, as well as related topics, rely crucially on the existence of particularly nice semisimple elements in the Lie algebra or in the group itself that were called ``good'' by Adler (\cite{Adler}), or on the corresponding notion of ``genericity'' introduced by Yu (\cite{Yu}) for the dual of the Lie algebra.
More precisely, a common assumption is that all tori of the reductive group split over a tamely ramified extension and that every Moy--Prasad filtration coset of (the Lie algebra of) such a torus contains a good element. However, it has been unclear under what conditions these assumptions are satisfied. In the present paper, we obtain sufficient conditions and necessary conditions on the reductive group $G$ and the residual characteristic $p$ for the assumption to be satisfied. In particular, we provide a criterion for the existence of good semisimple elements in all Moy--Prasad filtration cosets of the group and of its Lie algebra. While this question in the Lie algebra setting has previously been considered under more restrictive hypotheses by Adler and Roche (\cite{Adler-Roche}), to our knowledge the present paper presents the first result about the existence of good semisimple elements in all positive depth Moy--Prasad filtration cosets of the group.
Good elements in the Lie algebra and their analogue in the dual of the Lie algebra  %, which are called generic elements, 
are important ingredients in the constructions of supercuspidal representations by Adler and Yu. The goodness property is also used in the study of inclusion of Bruhat--Tits buildings of Levi subgroups, for the existence of good minimal K-types, and for the study of characters and orbital integrals as well as Hecke algebras, to name a few examples. 

To explain our results in more detail, let us assume in the introduction that $G$ is an absolutely simple reductive group over a non-archimedean local field $k$. General reductive groups are treated in the main part of the paper.

It is an easy exercise to show that if $p$ does not divide the order of the Weyl group $\Wabs$ of $G$, then all tori of $G$ split over a tamely ramified field extension of $k$, see Proposition \ref{Prop-tame-tori-easy}. However, it seemed unknown if this condition is also necessary. In Theorem \ref{Thm-tame-tori-less-easy} we show that this condition is necessary unless $G$ is a non-split inner form of type $A_n$ ($n \geq 2$), $D_l$ ($l$ a prime greater 4), or $E_6$. In these cases we provide an optimal upper bound for the set of primes $p$ for which all tori of $G$ are tamely ramified, and we prove that this upper bound is achieved for some inner form, see Theorem \ref{Thm-tame-tori-less-easy} Part \ref{Prop-Part-2}. 
This is done in Section \ref{section-tame-tori}.

In Section \ref{section-good-elements} we show that the condition $p \nmid \abs{\Wabs}$ is also sufficient for the existence of sufficiently many good semisimple elements in the following sense. % We call a semisimple element $X$ in the Lie algebra of $G$ a \textit{good} element if it is contained in the Lie algebra $\ft$ of a tame torus $T$ and all roots of $G$ with respect to $T$ that do not vanish on $X$ have the same valuation (at $X$). 
%In the 1990s, Moy and Prasad introduced filtrations of (the rational points of) reductive groups by compact open subgroups, which became an important tool in the study of representations of $p$-adic groups. 
If $T$ is a maximal torus of $G$ with Lie algebra $\ft$, then we denote its  Moy--Prasad filtration indexed by a discrete subset of the rational numbers $\hdots, r_{-1}, r_0, r_1, r_2, \hdots$ by
$$ \hdots \supsetneq \ft_{r_{-1}} \supsetneq  \ft_{r_0} \supsetneq \ft_{r_1} \supsetneq \ft_{r_2} \supsetneq \hdots .$$ 
Suppose $T$ splits over a tamely ramified extension of $k$. Then, for $n \in \bZ$, we call an element $X$ in $\ft_{r_n}$ a \textit{good} element if all roots of $G$ with respect to $T$ that do not vanish at $X$ have valuation $r_n$ (when evaluated on $X$).
The assumption frequently used in the representation theory of and harmonic analysis on the group $G(k)$ is that for every integer $n$, every coset in $\ft_{r_n}/\ft_{r_{n+1}}$ contains a good element. In Theorem \ref{Thm-good} we prove that this is always true under the assumption that $p \nmid \abs{\Wabs}$. In Theorem \ref{Thm-good-group} we prove the analogous statement in the group setting. 

As mentioned above, the Lie algebra setting has already been studied earlier by Adler and Roche under more restrictive hypotheses, see \cite[Section~5]{Adler-Roche}. For example, they assume that $ p > 113, 373$ or 1291 if $G$ is of type $E_6, E_7$ or $E_8$, respectively, while we only require $p > 5,  7$ or $7$ in these cases. However, their proof for general reductive groups (i.e. the proof of \cite[Lemma~5.12]{Adler-Roche}) seems to contain a mistake that led us to providing this new, more general proof.

%\textbf{things to be defined}:\\
%$k$ local non-arch field (of arbitrary char!)\\

%$k^s$ separable closure of $k$, and all separable extensions are viewed inside $k^s$\\
%$p$ residue field characteristic \\
%$q$ cardinality of the residue field of $k$\\
%$G$ (connected) reductive group over $k$\\
%$\Wabs$= Weyl group of $G_{k^s}$ \\
%$X_*(T_{k^s})$ = cocharacter group\\
%$X^*(T_{k^s})$\\
%$d \chi =\Lie(\chi)$\\
%$  N_G(T)$ \\
%$Z_G(S)$ \\
%$Z_G(S)^\circ$\\
%$\delta_{i,p}$= 1 or 0 -maybe no longer needed\\
% $\bQ\Phi$ - maybe no longer needed\\
% $Z(G)$\\
%$ \Cent_{\Wabs}(w_5)$\\

%$\Phi(G)$\\
%$\Phi(G,T_{\ov k})$\\
%$d \alpha$ (for $\alpha \in \Phi$)\\
%$\bZ \Phi_0 $ and $\bQ\Phi$ smallest ?-lattice inside ... containing $\Phi$\\
%$\val$ \\
%$\ft=\Lie(T)(k)$, $\ft_r$, $\ft(E)_{r}$\\
%$\check \Phi$\\
%$\<\check\alpha,\beta\>$, the pairing \\
% $s_\alpha$ element in weyl group for root alpha\\
% $H_{\alpha_1}$\\

%commonly used notation
%$\cO_E$, $\cO$\\
%$G_E, T_E$\\
%$\delta_{\alpha, \beta}$\\
%$\bZ_{(p)}$\\
%$S_3$  = the symmetric group on three letters \\
%$EF$ = composite field of $E$ and $F$  where $E$ and $F$ are Galois extensions of $k$, note that $EF$ is Galois over $k$\\

%Galois cohomology notation $H^1(k,W_T)$\\

\textbf{Notation.} Throughout the paper $k$ denotes a non-archimedean local field of arbitrary characteristic, and $k^*=k-\{0\}$. We denote by $k^s$ a separable closure of $k$, and all separable extensions of $k$ are viewed inside $k^s$. We let $\cP$ be the maximal ideal of the ring of integers $\cO$ of $k$, and $\cP_{k^s}$ the maximal ideal of the ring of integers of $k^s$. We write $\val: k \ra \bZ \cup \{\infty\}$ for the discrete valuation on $k$ with image $\bZ \cup \{\infty\}$, and we also denote by $\val: k^s \ra \bQ \cup \{\infty\}$ its extension to $k^s$. We denote the characteristic of the residue field $\cO/\cP$ by $p$ and its cardinality by $q$.

Unless mentioned otherwise, $G$ denotes a reductive group over $k$, and we use the convention that all reductive groups are connected.

If $S$ is a torus of $G$, we denote by $N_G(S)$ and $Z_G(S)$ the normalizer and centralizer of $S$ in $G$, respectively,  and we write $Z(G)$ for the center of $G$. If $T$ is a maximal torus of $G_{E}:=G \times_k E$ for some extension $E$ of $k$, then $\Phi(G,T)$ denotes the roots of $G_E$ with respect to $T$, and we might abbreviate $\Phi(G,T_{k^s})$ by $\Phi(G)$ if the choice of torus does not matter. We write $\Wabs$ for the Weyl group of $G(k^s)$, i.e. $\Wabs \simeq N_G(T(k^s))/T(k^s)$ for some maximal torus $T$ of $G_{k^s}$. If $T$ is a torus defined over some field $E$ (usually a finite extension of $k$, or $k^s$), then we write $X_*(T)=\Hom_E(\bG_m,T)$ for the group of cocharacters and $X^*(T)=\Hom_E(T,\bG_m)$ for the group of characters of $T$ defined over $E$. For a subset $\Phi$ of $X^*(T) \otimes_\bZ \bR$ (or $X_*(T)\otimes_\bZ \bR$) and $R$ a subring of $\bR$, we denote by $R\Phi$ the smallest $R$-submodule of $X^*(T) \otimes_\bZ \bR$ (or $X_*(T)\otimes_\bZ \bR$, respectively) that contains $\Phi$. For $\chi \in X^*(T)$, we denote by $d\chi \in \Hom_E(\Lie(T),\Lie(\bG_m))$ the induced morphism of Lie algebras. If $T$ is split, we use the notation $\< \cdot, \cdot\>: X_*(T) \times X^*(T) \ra \bZ$ for the standard pairing.
If $\Phi$ is a root system, then we denote by $\check \Phi$ the dual root system, and for $\alpha \in \Phi$, we let $\check \alpha$ be the corresponding dual root and $\check \omega_\alpha$ the fundamental coweight, i.e. $\check \omega_\alpha(\beta)=\delta_{\alpha,\beta}$ for $\beta \in \Phi$ (as usual, the expression $\delta_{a,b}$ denotes the number one if $a=b$ and zero otherwise.). Moreover, we denote by $s_\alpha$ the reflection corresponding to $\alpha$ in the Weyl group associated to $\Phi$, e.g. in $\Wabs$ if $\Phi=\Phi(G)$. If $T$ is a split maximal torus of $G_E$ defined over an extension $E$ of $k$, then for $\alpha \in \Phi(G,T)$, we set $H_\alpha=d\check\alpha(1) \in \Lie(T)(E)$.%, and $H_{\omega_\alpha}$ is defined to be the element in the $\bR$-subspace of $X_*(T)\otimes_{\bZ}\bR$ spanned by $\{H_\beta\, | \, \beta \in \Phi(G,T)\}$.

In addition, we use the following common notation: If $E$ is an extension of $k$, then $\cO_E$ denotes the ring of integers in $E$. As already used above, the base change of a scheme $S$ over $k$ to $E$ is denoted by $S_E$. If $E$ and $F$ are extensions of $k$, then we denote by $EF$ the composite field. We write $\bZ_{(p)}$ for the integers localized away from the prime ideal $(p)$, and $S_n$ for the symmetric group on $n$ letters for $n \in \bZ_{>0}$. If $w$ is an element of a finite group $W$, then $\Cent_W(w)$ denotes the centralizer of $w$ in $W$.
In addition, if $H$ is a group scheme over $k$, then we write $H^1(k,H)$ for the Galois cohomology set $H^1(\Gal(k^s/k),H(k^s))$, and $\Hom(A,B)/\mathord\sim$ denotes the set of group homomorphisms from a group $A$ to a group $B$ up to conjugation by elements of $B$.

\textbf{Acknowledgments.}  The author thanks Jeffrey Adler, Stephen DeBacker, Tasho Kaletha, Ju-Lee Kim and Alan Roche for discussions related to this paper and Beth Romano for feedback on an earlier version of this paper.
The author also thanks the Max-Planck-Institut für Mathematik and the Institute for Advanced Study for their hospitality and wonderful research environment.  While at the Institute for Advanced Study the author was supported by the National Science Foundation under Grant No. DMS - 1638352. The author also appreciates earlier support by a postdoctoral fellowship of the German Academic Exchange Service (DAAD).

\section{Condition for tame tori}
\label{section-tame-tori}
%In this section we will obtain necessary and sufficient conditions to guarantee that all tori of $G$ split over a tamely ramified extension of $k$. 
We begin with the case that $G$ is an absolutely simple reductive group, and recall the following most likely well known result.
\begin{Prop} \label{Prop-tame-tori-easy}
	Let $G$ be an absolutely simple reductive group defined over a non-archimedean local field $k$ of residual characteristic $p$. If $p$ does not divide $\abs{\Wabs}$, the order of the Weyl group $\Wabs$ of $G(k^s)$, then every torus of $G$ splits over a (finite) tamely ramified field extension of $k$.
\end{Prop}
\Proof
Let $T$ be a torus of $G$, and assume $p \nmid \abs{\Wabs}$. % and show that $T$ splits over a tamely ramified extension. 
Since every torus is contained in a maximal torus that is also a maximal torus over $k^s$,
%By applying Grothendieck's theorem (\coR{\cite[XIV, 1.1]{SGA3}}) on the existence of tori over $k$ that are maximal over $k^s$ to the centralizer of $T$, 
 we may assume that $T_{k^s}$ is a maximal torus of $G_{k^s}$. Then the $\Gal(k^s/k)$-action on $T_{k^s}$ factors through $\Wabs \rtimes \Aut(\Dyn)$, where $\Aut(\Dyn)$ denotes the automorphisms of the Dynkin diagram corresponding to $\Phi(G)$. Thus, since $p \nmid \abs{\Wabs}$ and $\Aut(\Dyn) \mid \abs{\Wabs}$, hence $p \nmid \abs{\Wabs \rtimes \Aut(\Dyn)}$, the $\Gal(k^s/k)$-action factors through the action of $\Gal(E/k)$ for some finite tamely ramified field extension $E$ of $k$, and $T$ splits over $E$. \qed

In order to address the reverse direction, which turns out to be true as well except for some non-split inner forms of split groups of type $A_n, D_n$ and $E_6$, we begin with two lemmata. The first lemma shows the existence of certain Galois extensions of $k$ that will be used to construct tori that do no split over any tamely ramified extension of $k$ when certain conditions are met.
\begin{Lemma} \label{Lemma-deg-p-extension}
	Let $k$ be a non-archimedean local field of residual characteristic $p$.
\begin{enumerate}[label=(\arabic*),ref=\arabic*]
	\item \label{item-extensions-0} Let $N$ be an integer divisible by $p$. Then there exists a wildly ramified Galois extension $E$ of $k$ with Galois group $\bZ/N\bZ$. 
%	\item\label{item-extensions-1}	There exists a (totally, wildly) ramified Galois extension $E$ of $k$ of degree $p$.
%	\item  \label{item-extensions-2} There exists a (wildly) ramified Galois extension $E$ of $k$ with Galois group $\Gal(E/k)\simeq \bZ/p^2 \bZ$.
	\item \label{item-extensions-3} Suppose $F$ is a finite Galois extension of $k$ with Galois group $\bZ/1\bZ, \bZ/2\bZ, \bZ/3\bZ$ or $S_3$. Then there exists a Galois extension $E$ of $k$ of degree 2 such that $E \cap F=k$, hence $\Gal(EF/k)=\Gal(E/k)\times\Gal(F/k)$.
\end{enumerate}
\end{Lemma}
\Proof By local class field theory Galois extensions of $k$ with cyclic Galois group $\Gamma$ correspond to open subgroups $U$ of $k^*$ such that $k^*/U \simeq \Gamma$. Since $k$ has precisely one unramified extension of degree $N$, it suffices to observe that $k^*$ has at least two open subgroups $U$ such that $k^*/U \simeq \bZ/N\bZ$ to prove \eqref{item-extensions-0}. This follows from the following isomorphisms of topological groups (\cite[Chapter~II.~(5.7.)~Proposition]{Neukirch}): If $k$ has characteristic zero, then
$$ k^* \simeq \bZ \oplus \bZ/(q-1)\bZ \oplus \bZ/p^a\bZ \oplus \bZ_p^d$$
for some $a \geq 0$ and $d=[k:\bQ_p]$, and if $k$ has characteristic $p$, then
$$ k^* \simeq \bZ \oplus \bZ/(q-1)\bZ\oplus \bZ_p^\bN .$$
More precisely, the unramified extension of degree $N$ corresponds to the subgroup $N\bZ$ in the first summand, and an example of a wildly ramified Galois extension of $k$ with Galois group $\bZ/N\bZ$ is provided by taking the one corresponding to the subgroup given by $N_{p'}\bZ$ in the first summand, $N_{p}\bZ_p$ in the last summand and all the full remaining summands, where $N=N_{p'}N_p$ with $p \nmid N_{p'}$ and $N_{p}$ is a power of $p$.

In order to prove \eqref{item-extensions-3}, first consider the case that $F$ is either unramified or totally ramified. In this case we can take $E$ to be a degree two totally ramified or unramified extension of $k$, respectively, which exists by considerations as in the proof of \eqref{item-extensions-0} (note that $q-1$ is even if $p \neq 2$). Then $E \cap F =k$, as desired. If $F$ is neither unramified nor totally ramified, then $\Gal(E/k) \simeq S_3$. The group $S_3$ contains a unique order three normal subgroup, which corresponds to a degree two Galois extension $E'$ of $k$ contained in $E$.  If $E'$ is unramified, then we choose $F$ to be a degree two totally ramified Galois extension, and otherwise we let $F$ be the degree two unramified extension of $k$. In both cases $E \cap F = E' \cap F = k$.
\qed

For the second lemma, we allow $G$ to be any (connected) reductive group over $k$. We fix a maximal torus $T$ of $G$ and let $W_T=N_G(T)/T$. If $T'$ is any torus of $G$, then there exists $g \in G(k^s)$ such that $T'=gTg^{-1}$, and the map $\Gal(k^s/k) \ra N_G(T)(k^s)$ given by $\sigma \mapsto g^{-1}\sigma(g)$ defines a cocycle whose image in $H^1(k,N_G(T))$ is independent of the choice of $g$. On the other hand, any cocycle $\sigma \mapsto n_\sigma$ in $H^1(k,N_G(T))$ that maps to the trivial cocycle in $H^1(k,G)$ gives rise to a (conjugacy class of a) torus $T'$ in $G$ as follows. Write $n_\sigma=g^{-1}\sigma(g)$ for some $g \in G(k^s)$ and set $T'=gTg^{-1}$. Then for all $\sigma \in \Gal(k^s/k)$ we have $$\sigma(T')=\sigma(g)\sigma(T)\sigma(g^{-1})=gn_\sigma T n_\sigma^{-1} g^{-1}=gTg^{-1}=T', $$
hence the tours $T'$ is defined over $k$. Moreover, a different choice $g'\in G(k^s)$ with $n_\sigma=(g')^{-1}\sigma(g')$ yields a torus that is $G(k)$-conjugate to $T'$. Raghunathan proved in \cite[Main~Theorem~1.1]{Raghunathan} that the map 
$$ H^1(k,N_G(T)) \ra H^1(k,W_T)$$
induced by the quotient map $N_G(T) \ra W_T=N_G(T)/T$ is surjective and that if $G$ is quasi-split and $T$ is the centralizer of a maximal split torus, then every element of $H^1(k,W_T)$ lifts to an element in $\ker(H^1(k,N_G(T)) \ra H^1(k,G))$. A similar result was obtained around the same time by Gille (\cite{Gille}). 
As a consequence we obtain the following lemma, for which an alternative proof is given by Kaletha in  \cite[Lemma~3.2.2]{Kaletha}.
\begin{Lemma} \label{Lemma-transfer-to-quasi-split}
	Let $G$ be a (connected) reductive group defined over a non-archimedean local field $k$, and let $G'$ be its quasi-split inner form. If $T$ is a maximal torus of $G$, then $T$ is isomorphic to a maximal torus of $G'$.  
\end{Lemma}
\Proof\\
Let $\pi: G \ra G^{ad}=G/Z(G)$ be the adjoint quotient of $G$. Then the maximal tori of $G$ are in one to one correspondence with maximal tori in $G^{ad}$ via $T \mapsto \pi(T)$ (with inverse $T \mapsto \pi^{-1}(T)$). Since the action of $\Gal(k^s/k)$ on tori is determined the action on the center of $G$, which is the same for all inner forms, and by the permutation of root groups, it suffices to treat the case of $G$ being adjoint.

Let $T$ be a maximal torus of $G$ and $T'$ the centralizer of a maximal split torus of $G'$. Let $\phi:G'_{k^s}\ra G_{k^s}$ be an isomorphic that sends $T'$ to $T$. Then for $\sigma \in \Gal(k^s/k)$, the automorphism $\phi^{-1}\circ \sigma(\phi)$ corresponds to conjugation by an element $n_\sigma \in N_G(T)(k^s)$, which yields a cocycle $\sigma \mapsto n_\sigma$ in $H^1(k,N_G(T))$. By \cite[Main~Theorem~1.1]{Raghunathan} there exists a cocycle $\sigma \mapsto n'_\sigma \in H^1(k,N_G(T))$ whose image in $H^1(k,G)$ is trivial and whose image in $H^1(k,W)$ agrees with the image of $\sigma \mapsto n_\sigma$. Let $g \in G(k^s)$ such that $n'_\sigma =g^{-1}\sigma(g)$ for all $\sigma \in \Gal(k^s/k)$. Then $T''=gT'g^{-1}$ is a maximal torus of $G'$ that is defined over $k$ and is isomorphic to $T$, because the image of $\sigma \ra n_\sigma$ and $\sigma \ra n_\sigma'$ under $H^1(k,N_G(T')) \ra H^1(k,W_{T'}) \ra H^1(k,W_{T'}) \ra H^1(k,\Aut(T'))$ agree.
\qed

Now we are ready to prove the reverse direction to Proposition \ref{Prop-tame-tori-easy} and characterize the set of primes $p$ for which $G$ contains only tame tori.

\begin{Thm} \label{Thm-tame-tori-less-easy}
	Let $G$ be an absolutely simple reductive group defined over a non-archimedean local field $k$ of residual characteristic $p$. 
\begin{enumerate}
		\item \label{Prop-Part-1}
If $G$ is not of one of the following forms
\begin{itemize}
	\item a non-split inner form of a split group of type $A_n$ for some positive integer $n \geq 2$,
	\item a non-split inner form of a split group of type $D_l$ for some prime number $l \geq 4$,
	\item a non-split inner form of a split group of type $E_6$,
\end{itemize}		
	then the following two statements are equivalent 
\begin{enumerate}[label=(\roman*),ref=\roman*]
	\item \label{i} Every torus of $G$ splits over a (finite) tamely ramified field extension of $k$.
	\item \label{ii} $p$ does not divide $\abs{\Wabs}$.
\end{enumerate}	

\item \label{Prop-Part-2} If $G$ is 
	\begin{itemize}
		\item an inner form of a split group of type $A_n$ for some positive integer $n \geq 2$, or
		\item an inner form of a split group of type $D_l$ for some prime number $l \geq 4$, or
		\item an inner form of a split group of type $E_6$,
	\end{itemize}	
	then the set of primes $\cS_G$ for which every torus of $G$ splits over a tamely ramified field extension of $k$ is not the same for all inner forms of $G$. 
	
	More precisely, denote by $\cS_G^{-}$ the set of all primes $p$ with $p \nmid \abs{\Wabs}$, and 
\begin{enumerate}[ref=\alph*]	
	\item \label{item-A-n} if $G$ is an inner form of a split group of type $A_n$, then let $\cS_G^{+}$ be the set of all primes $p$ that do not divide $n+1$,

	\item \label{item-D-l} if $G$ is an inner form of a split group of type $D_l$, then let $\cS_G^{+}$ be   the set of all primes $p$ that do not divide  $\abs{\Wabs}/l=2^{l-1}(l-1)!$,

	\item \label{item-E-6} if $G$ is an inner form of a split group of type $E_6$, then let $\cS_G^{+}$ be  the set of all primes $p$ that do not divide  $\abs{\Wabs}/5=2^{7}\cdot3^4$.
\end{enumerate}	
	Then we have $\cS_G^{-} \subseteq \cS_G \subseteq \cS_G^{+}$, and there exist inner forms $G^-$ and $G^+$ of $G$ such that $\cS_{G}^-=\cS_{G^-}$ and $\cS_{G}^+=\cS_{G^+}$.
\end{enumerate}	
\end{Thm}

\begin{Rem} \label{Remark-inner-forms}
 In the case of $G$ being an inner form of a split group of type $A_n$, $D_l$ or $E_6$, we do not state which set of primes has to be excluded for which inner forms. However, the interested reader might use the methods employed in the following proof to work out a classification.	
\end{Rem}

For the readers convenience, we list the order of the irreducible Weyl groups in Table \ref{table-weyl-group} (\cite[VI.4.5-VI.4.13]{Bourbaki-4-6}).

\begin{table}[h]\footnotesize
\begin{tabular}{|c|c|c|c|c|c|c|c|c|}
	\hline type  & 			$A_n \, (n \geq 1)$     & $B_n, C_n \, (n \geq 2)$  &  $D_n \, (n \geq 3)$              & $E_6$                 &  $E_7$                      &  $E_8$ & $F_4 $ & $G_2$  \\ 
	\hline $\abs{\Wabs}$ &  $(n+1)!$  & $2^n \cdot n!$   &  $2^{n-1} \cdot n!$ & $2^7\cdot3^4\cdot5$  &  $2^{10}\cdot3^4\cdot5\cdot7$ & $2^{14}\cdot3^5\cdot5^2 \cdot 7$ &  $2^7 \cdot 3^2$ & $2^2 \cdot 3$ \\ 
	\hline 
\end{tabular} 
	\caption{Order of irreducible Weyl groups (\cite[VI.4.5-VI.4.13]{Bourbaki-4-6})}
	\label{table-weyl-group}
\end{table}

\textbf{Proof of Theorem \ref{Thm-tame-tori-less-easy}.}\\
\textbf{Proof or Part \ref{Prop-Part-1}.} By Proposition \ref{Prop-tame-tori-easy} it suffices to show that if \eqref{ii}  is not satisfied, than \eqref{i} is not true either, so let us assume that $p$ divides $\abs{\Wabs}$. Recall that over a non-archimedean local field every anisotropic maximal torus transfers to all inner forms (see \cite[\S10]{Kottwitz} and \cite[Lemma~3.2.1]{Kaletha}). Hence it suffices to exhibit an anisotropic maximal torus that does not split over a tamely ramified extension for the case when $G$ is quasi-split, so we assume that $G$ is quasi-split for the remainder of the proof of Part \ref{Prop-Part-1}. Let $T$ be a maximal torus of $G$ contained in a Borel subgroup $B$ of $G$ (that is defined over $k$). Let $W=N_G(T)/T$, hence $\Wabs=W(k^s)$.
We distinguish three cases.

\textit{Case 1: $\Wabs$ is not of type $A_n, D_{2n+1}$ or $E_6$ with $n \in \bZ_{\geq 2}$.} In other words, $\Wabs$ contains the element -1 when viewed as a reflection subgroup acting on $X_*(T_{k^s}) \otimes_\bZ \bR$ for some maximal torus $T$ of $G$. The equivalence of these two conditions follows from \cite[3.19~Corollary]{Humphreys} combined with  \cite[Table~3.1]{Humphreys}.

 If $T$ is split, then $\Wabs=W(k)$. If $T$ is not split, then $\Wabs$ is of type $D_{2n}$ for some $n\geq 2$, and $W(k)=\Wabs^{\Gal(k^s/k)}$, where the $\Gal(k^s/k)$-action factors through the group of automorphism of $\Wabs$ induced by a Dynkin diagram automorphisms. If the $\Gal(k^s/k)$-action factors through an order two automorphism of $\Wabs$, then, by \cite[Theorem~32 and p.~175]{SteinbergYale}, the group $W(k)$ is a Weyl group of type $B_{2n-1}$. Thus, if $p \mid \abs{\Wabs}=2^{n-1}(2n)!$, then $p \mid \abs{W(k)}=2^{n-1}(2n-1)!$. If the $\Gal(k^s/k)$-action on $\Wabs$ does not factor through an order two automorphism, then $n=2$, and it has to factor through the action of $S_3$, and the normal subgroup $A_3$ of order three has to act non-trivially. By \cite[Theorem~32 and p.~176]{SteinbergYale}, the fixed subgroup of $\Wabs$ under the action of $A_3$ is a Weyl group of type $G_2$, which is generated by $s_3$ and $s_1s_2s_4$, where $s_1, s_2, s_3, s_4$ denotes a set of simple reflections of $\Wabs$ (determined by a Borel subgroup containing $T$) such that $A_3$ permutes the commuting elements $s_1, s_2, s_4$. Since $S_3$ acts by permuting $\{s_1, s_2, s_4 \}$ and fixing $s_3$, we have $\Wabs^{A_3}=\Wabs^{S_3}=W(k)$. Hence, $p \mid \abs{\Wabs}=2^3\cdot 4!$ implies that $p \mid \abs{W(k)}=2^2\cdot3$. Thus, in all cases $p \mid \abs{W(k)}$, and we choose an element $w$ of order $p$ in $W(k)$.

In addition, we have $-1 \in W(k)=\Wabs^{\Gal(k^s/k)}$, because the Galois action preserves the center of $\Wabs$, which only consists of 1 and -1 (\cite[6.3~Proposition~(d)]{Humphreys}).

Let $E_1$ be a minimal Galois extension of $k$ over which $T$ splits. We assume that $E_1$ is tamely ramified over $k$ because otherwise we are done. Note that we observed above that $\Gal(E_1/k)$ is isomorphic to $\bZ/1\bZ$, $\bZ/2\bZ$, $\bZ/3\bZ$ or $S_3$. Let $E_2$ be a degree two Galois extension of $k$ such that $E_1 \cap E_2 =k$, which exists by Lemma \ref{Lemma-deg-p-extension}\eqref{item-extensions-3}.  If $E_2$ is tamely ramified over $k$, then let $E_3$ be a totally, wildly ramified Galois extension of $k$ of degree $p$, which exists by Lemma \ref{Lemma-deg-p-extension}\eqref{item-extensions-0}, otherwise let $E_3=k$. Then we have $E_1E_2 \cap E_3 =k$, and hence $\Gal(E_1E_2E_3/k) \simeq \Gal(E_1/k) \times \Gal(E_2/k) \times \Gal(E_3/k)$. Note that $\Gal(E_2/k) \times \Gal(E_3/k)$ acts trivially on $\Wabs$ and $\Gal(k^s/k)$ acts trivially on $W(k)=\Wabs^{\Gal(k^s/k)}$. This allows us to define a cocycle $f: \Gal(k^s/k) \ra \Wabs$ by requiring that it factors through $\Gal(E_2E_3/k)\simeq \bZ/2\bZ \times \bZ/[E_3:k]\bZ$ and satisfies $f((a,b))=(-1)^a \cdot w^{\frac{p}{[E_3:k]}\cdot b}$ for $(a,b) \in \bZ/2\bZ \times \bZ/[E_3:k]\bZ \simeq \Gal(E_2E_3/k)$. This is a well defined cocycle by the above discussion and it has nontrivial image in $H^1(k,W)$, which lifts to a cocycle $\sigma \mapsto n_\sigma$ in $H^1(k,N_G(T))$ that maps to the trivial class in $H^1(k,G)$ by \cite[Main~Theorem~1.1]{Raghunathan}. Thus there exists $g \in G(k^s)$ such that  $n_\sigma=g^{-1}\sigma(g)$ for all $\sigma \in \Gal(k^s/k)$, and $T'=gTg^{-1}$ is defined over $k$. Moreover, since the nontrivial element of $\Gal(E_2/k)$ acts via -1 on $X_*(T'_{k^s})$, the maximal torus $T'$ is anisotropic. It remains to show that $T'$ does not split over a tamely ramified extension. Suppose $T'$ split over a tamely ramified extension $F$, and let $E$ be the composite $E_1F$. Then the image of $\sigma \mapsto n_\sigma$ in $H^1(E,N_G(T))$ would be the trivial class, hence (the image of) $f$ would represent the trivial class in $H^1(E,W)=\Hom(\Gal(k^s/E),W(E))/\mathord\sim$. This would contradict the fact that the image of $\Gal(k^s/E)$ under (the image of) $f$ in $\Hom(\Gal(k^s/E),W(E))$ contains the nontrivial element $w$ (if $E_3 \neq k$) or -1 (if $E_3=k$).  
\\

\textit{Case 2: $\Wabs$ is of type $A_n, D_{2n+1}$ or $E_6$ with $n \in \bZ_{\geq 2}$ and $G$ is not split.} We assume $T$ splits over a tamely ramified extension as otherwise we are done. Then the action of $\Gal(k^s/k)$ on $T$ factors through a tamely ramified quadratic Galois extension $E_1/k$. We denote the induced action of the non-trivial element $\rho \in \Gal(E_1/k)$ on $X_*(T_{k^s})$ by $\tau$. Then $\tau$ is the automorphism induced by the non-trivial Dynkin-diagram automorphism (for the set of simple roots $\Delta$ determined by the Borel subgroup $B$). Let $w_0$ be the longest element in $\Wabs$ (with respect to  $\Delta$). Then $w_0(\Delta)=-\Delta$ and $w_0^2=1$.  %, and hence the action of $w_0$ on $\Wabs$ via conjugation is an order two automorphism that we denote by $\Inn(w_0)$.
 For $\alpha \in \Delta$, denote by $s_\alpha$ the corresponding simple reflection in $\Wabs$. Then $w_0 s_\alpha w_0=s_{-w_0(\alpha)}$, and hence the automorphism of $\Wabs$ obtained by conjugation by $w_0$ stabilizes the set of simple reflections of $\Wabs$, i.e. arises from an automorphism of the Dynkin diagram. Since $-1 \not \in \Wabs$, we have $w_0 \neq -1$, and therefore $w_0 = - \tau$.
 
Let $E_3$ be a totally, wildly ramified Galois extension of $k$ of degree $p$, which exists by Lemma \ref{Lemma-deg-p-extension}\eqref{item-extensions-0}, and let $w \in \Wabs$ be an element of order $p$. We define $f: \Gal(k^s/k) \ra \Wabs$ by requiring that it factors through $\Gal(E_1E_3/k)\simeq \bZ/2\bZ \times \bZ/p\bZ$ and satisfies $f((a,b))=w^b \cdot w_0^a$ for $(a,b) \in \bZ/2\bZ \times \bZ/p\bZ \simeq \Gal(E_1E_3/k)$. Since $\Gal(k^s/E_1E_3)$ acts trivially on $\Wabs$, and 
$$ f((a,b)+(c,d)) =  w^{b+d} \cdot w_0^{a+c} =    w^b \cdot w_0^a \cdot (w_0^a (w^d \cdot w_0^c) w_0^a) = f((a,b)) \cdot \rho^a(f((c,d)))$$
for $(a,b)$ and $(c,d)$ in $\bZ/2\bZ \times \bZ/p\bZ \simeq \< \rho \> \times \bZ/p\bZ  \simeq \Gal(E_1E_3/k)$, the map $f$ is a cocycle that defines an element in $H^1(k,W)$. As in Case 1, by \cite[Main Theorem~1.1]{Raghunathan} this element of $H^1(k,W)$ lifts to a cocycle $\sigma \mapsto g^{-1}\sigma(g)$ in $H^1(k, N_G(T))$ for some $g \in G(k^s)$, and $T'=gTg^{-1}$ is a torus defined over $k$. Let $\rho' \in \Gal(k^s/k)$ be a preimage of $(\rho,0) \in \Gal(E_1/k) \times \bZ/p\bZ \simeq \Gal(E_1/k) \times \Gal(E_3/k)$. Then $\rho'$ acts on $X_*(T'_{k^s})$ via $w_0 \tau=w_0 (-w_0)=-1$, and hence $T'$ is anisotropic. Moreover, $T'$ does not split over a tamely ramified extension by construction (and the same argument as in Case 1), and therefore \eqref{i} does not hold, which we needed to show.
\\

\textit{Case 3: $\Wabs$ is of type $A_n, D_{2n+1}$ or $E_6$ with $n \in \bZ_{\geq2}$ and $G$ is split.} If $2n+1$ is not a prime number, then all primes $p$ dividing $\abs{\Wabs}=2^{2n}(2n+1)!$ for $\Wabs$ of type $D_{2n+1}$ also divide $\abs{\Wabs}/(2n+1)=2^{2n}(2n)!$, and hence the first paragraph of the proof of Part  \ref{Prop-Part-2}\eqref{item-D-l} (which applies to $l=2n+1$ not being a prime number as well) together with Proposition \ref{Prop-tame-tori-easy} will show that \eqref{i} and \eqref{ii} are equivalent. Hence, by the theorem statement, the remaining cases to treat are those for which $G$ splits. Therefore it suffices to exhibit a torus $T' \subset G$ that does not split over a tamely ramified extension. Let $E$ be a degree $p$ totally ramified Galois extension  of $k$, and let $w \in \Wabs$ be an element of order $p$. Then we can define a cocycle $f:\Gal(k^s/k) \ra \Wabs$ by requiring that it factors through $\Gal(E/k)$ and sends a generator of $\Gal(E/k)$ to $w$. This yields an element of $H^1(k,W)$, which can be lifted to a cocycle $\sigma \mapsto g^{-1}\sigma(g) \in H^1(k,N_G(T))$ for some $g \in G(k^s)$ by \cite[Main~Theorem~1.1]{Raghunathan}, and $T'=gTg^{-1}$ is a desired torus (by the same reasoning as in Case 1).\\

%This concludes the proof of Part \ref{Prop-Part-1}.

\textbf{Proof of Part \ref{Prop-Part-2}.}\\
Note that $\cS_G^- \subseteq \cS_G$  by Proposition \ref{Prop-tame-tori-easy}, and if $G^-$ denotes the split inner form of $G$, then $\cS_{G}^-=\cS_{G^-}$  by Part \ref{Prop-Part-1}. In order to show that $\cS_G \subseteq \cS_G^+$ and that $\cS_{G}^+=\cS_{G^+}$ for some inner form $G^+$ of $G$, we treat the three cases \eqref{item-A-n}, \eqref{item-D-l} and \eqref{item-E-6} separately.

\textit{\eqref{item-A-n}} 
Let $G$ be an inner form of a split group of type $A_n$ for some positive integer $n \geq 2$, and let $\pi: G^{sc} \ra G$ be the map from the simply connected cover of $G$ to $G$. Then the maximal tori of $G^{sc}$ are in one to one correspondence with the maximal tori in $G$ via $T \mapsto \pi(T)$. Since the action of $\Gal(k^s/k)$ on tori is determined by the permutation of root groups, the splitting fields of $T$ and $\pi(T)$ coincide. Therefore it suffices to treat the case when $G$ is simply connected, i.e. $G$ is an inner form of $\SL_{n+1}$.

The anisotropic maximal tori of $\SL_{n+1}$ are isomorphic to the norm-one subtori of $\Res_{E/k}\bG_m$ for a degree-$(n+1)$ separable extensions $E/k$. Hence, if $p \mid n+1$, then we can choose a degree-$(n+1)$ wildly ramified Galois extension $E$ of $k$ by Lemma \ref{Lemma-deg-p-extension}\eqref{item-extensions-0}. Then the norm-one subtorus of $\Res_{E/k}\bG_m$ is an anisotropic maximal torus of $\SL_{n+1}$, hence transfers to $G$, and does not split over a tame extension of $k$. Thus, if every torus of $G$ splits over a tamely ramified field extension of $k$, then $p \nmid n+1$. 
 
 On the other hand, consider the case that $G$ is an anisotropic inner form of $\SL_{n+1}$. Then all tori in $G$ are anisotropic, and by Lemma \ref{Lemma-transfer-to-quasi-split} isomorphic to tori of $\SL_{n+1}$.
In addition, if $p \nmid n+1$, then any degree-$(n+1)$ separable extension of $k$ is tame and contained in a tamely ramified Galois extension of $k$ over which $\Res_{E/k}\bG_m$ (and its norm-one subtorus) splits. Hence all tori of $G$ split over a tamely ramified extension of $k$ if $p \nmid n+1$. 
 \\
 
\textit{\eqref{item-D-l}} Let $G$ be an inner form of a split group of type $D_l$ for some prime number $l \geq 4$. By \cite[Proposition~25]{Carter} together with \cite[Table~3]{Carter} for $1 \leq i \leq \frac{l-1}{2}$ there exists an element $w_i$ in $\Wabs$ whose characteristic polynomial is given by $(t^{l-i}+1)(t^i+1)$ (when viewing $\Wabs$ as acting on $X_*(T_{k^s}) \otimes_\bZ \bR$ for some maximal torus $T$ of $G$). Thus $w_i$ is elliptic, i.e.  $(X_*(T_{k^s}) \otimes_\bZ \bR)^{w_i}=\{0\}$, and the order $N_i$ of $w_i$ is the least common multiple of $2(l-i)$ and $2i$. Thus for $p \mid \abs{\Wabs}/l=2^{l-1}(l-1)!$, i.e. $p < l$, we can choose $1 \leq i \leq \frac{l-1}{2}$ so that $p \mid N_i$. Let $E$ be a wildly ramified Galois extension of $k$ with Galois group $\bZ/N_i\bZ$, which exists by Lemma \ref{Lemma-deg-p-extension} \eqref{item-extensions-0}. Then we can define a cocycle $f:\Gal(k^s/k) \ra \Wabs$ by requiring that $f$ factors through $\Gal(E/k)$ and sends a generator of $\Gal(E/k)$ to $w_i$. Analogous to Part \ref{Prop-Part-1} this yields an anisotropic torus in the split inner form of $G$  that therefore transfers to $G$ and that does not split over a tamely ramified extension.

It remains to show that for $p=l$ there exists an inner form $G$ of every split group of type $D_l$ for which all tori split over a tamely ramified extension. We consider the inner form $G$ corresponding to the symbol $^1D^{(2)}_{p,(p-3)/2}$ using the notation of \cite{Tits-index}. This group has the name $^4D_p$ in \cite{Tits}, and it has split rank $\frac{p-3}{2}$. Let $S$ be a maximal split torus of $G$. From the index in Figure \ref{Figure-index-D} provided by \cite[page~56]{Tits-index}, we see that the root system of $Z_G(S) \times_k k^s$ is of type $A_1^{(p-3)/2} \times A_3$.

\begin{figure}[h]
	\centering
	\begin{tikzpicture}
	
	\draw (0,0) -- (3.7,0);
	\draw (4.3, 0) -- (7,0);
	\draw (7,0) -- (8,-.5);
	\draw (7,0) -- (8,.5);
	
	\draw[fill=black] (0,0) circle(.1);% node [above] {$\alpha_1$};
	\draw[fill=black] (1,0) circle(.1);% node [above] {$\alpha_2$};
	\draw (1,0) circle(.2);
	\draw[fill=black] (2,0) circle(.1);% node [above] {$\alpha_3$};
	\draw[fill=black] (3,0) circle(.1);% node [above] {$\alpha_4$};
	\draw (3,0) circle(0.2);
	\draw (4,0)  node {$\cdots$};
	\draw[fill=black] (5,0) circle(.1);% node [above] {$\alpha_{p-4}$};
	\draw[fill=black] (6,0) circle(.1);% node [above] {$\alpha_{p-3}$};
	\draw (6,0) circle(.2);
	\draw[fill=black] (7,0) circle(.1);% node [above] {$\alpha_{p-2}$};
	\draw[fill=black] (8,-.5) circle(.1);% node [above] {$\alpha_{p}$};
	\draw[fill=black] (8,.5) circle(.1);% node [above] {$\alpha_{p-1}$};
	
	%	\node at (3,0) {$\cdots$};
	%	\node at (1.5, .3) {4};
	\node at (0,.4) {$\alpha_1$};
	\node at (1,.4) {$\alpha_2$};
	\node at (2,.4) {$\alpha_3$};
	\node at (3,.4) {$\alpha_4$};	
	\node at (5,.4) {$\alpha_{p-4}$};
	\node at (6,.4) {$\alpha_{p-3}$};
	\node at (7,.4) {$\alpha_{p-2}$};
	\node at (8,.8) {$\alpha_{p-1}$};
	\node at (8,-.2) {$\alpha_p$};
	
	\end{tikzpicture}
	\caption{Index for $^1D^{(2)}_{p,\frac{p-3}{2}}$ (\cite[page~56]{Tits-index})}	
	\label{Figure-index-D}	
\end{figure}

Suppose there exists a maximal torus $T$ of $G$ that does not split over a tamely ramified extension of $k$. By Lemma \ref{Lemma-transfer-to-quasi-split} we can identify $T$ with a maximal torus of the split inner form $G'$ of $G$, which yields a cocycle $f$ in $H^1(k,W_{T'}) \simeq \Hom(\Gal(k^s/k),\Wabs)/\mathord\sim$, where $W_{T'}=N_{G'}(T')/T'$ for some split maximal torus $T'$ of $G'$. Let $E$ be a Galois extension of $k$ such that $\Gal(k^s/E)$ is the kernel of (a representative of) $f \in \Hom(\Gal(k^s/k),\Wabs)$, and let $E^t$ be the maximal tamely ramified extension of $k$ contained in $E$. Then $\Gal(E/k)$ is isomorphic to a subgroup of $\Wabs$, and hence the normal $p$-subgroup $\Gal(E/E^t)$ is isomorphic to $\bZ/p\bZ$.
Without loss of generality (since all order-$p$ subgroups of $\Wabs$ are conjugate), we may assume that the image of $\Gal(E/E^t)$ is generated by $w=s_{\alpha_1}s_{\alpha_2}\hdots s_{\alpha_{p-1}}$ for simple roots $\alpha_1, \alpha_2, \hdots, \alpha_p$  as in Figure \ref{Figure-index-D}. This means $w$ is a Coxeter element in the subgroup $S_{p} \triangleleft \Wabs$ of type $A_{p-1}$ generated by $s_{\alpha_1}, \hdots, s_{\alpha_{p-1}}$. In other words, $w$ corresponds to a cycle of length $p$ in $S_p$. Suppose $w'$ is another element in the image of $\Gal(E/k)$ in $\Wabs$. Then $w'$ normalizes $\<w\>$, i.e. conjugation by $w'$ sends the cycle $w$ of length $p$ in $S_p$ to another cycle of length $p$ in $S_p$. Hence there exists $w'' \in S_p \triangleleft \Wabs$ such that $w'w''$ centralizes $w$. Thus the element $w'w''w$ has order a multiple of $p$. However, $\Wabs$ does not contain any elements of oder $Np$ for $N>1$ (see, e.g., \cite[Proposition~25]{Carter}). Thus $w'w'' \in \<w\>$, hence $w'\in S_p$, and the image of $\Gal(k^s/k)$ in $\Wabs$ is contained in $S_p=\< s_{\alpha_1}, \hdots, s_{\alpha_{p-1}} \>$.
Note that since $w$ is a Coxeter element of $S_p$, it has only 0 as fixed point when restricted to the subspace of $X_*(T_{k^s}) \otimes \bR$ spanned by $\check\alpha_1, \hdots, \check\alpha_{p-1}$. Hence $(X_*(T_{k^s})\otimes \bR)^{\Gal(k^s/k)}=\bR \cdot \check\omega_{\alpha_p}$, and $T$ has split rank one. If we denote by $S_T$ the maximal split subtorus of $T$, then by \cite[2.14~Lemma]{Steinberg-torsion} the group $Z_G(S_T)\times_k k^s$ is a reductive group whose root system consists of all the roots of $\Phi(G, T_{k^s})$ that are a linear combination of $\alpha_1, \alpha_2, \hdots, \alpha_{p-1}$, hence it is of type $A_{p-1}$. On the other hand, all maximal split tori of $G$ are conjugated over $k$, i.e. without loss of generality $S_T$ is contained in $S$. This implies that the split group $Z_G(S_T)\times_k k^s$ of type $A_{p-1}$ contains the split group $Z_G(S)\times_k k^s$ of type $A_1^{(p-3)/2}\times A_3$ as a Levi subgroup, which is not possible looking at their Dynkin diagrams. Hence all tori of $G$ split over a tamely ramified extension of $k$.
\\

\textit{\eqref{item-E-6}} Let $G$ be an inner form of a split group of type $E_6$. We proceed in the same spirit as in \eqref{item-D-l}. Let $w \in \Wabs$ be a Coxeter element. Then $w$ has order 12, and $(X_*(T_{k^s})\otimes \bR)^w=\{0\}$. Thus, if $p \mid \frac{\abs{\Wabs}}{5}=2^7 \cdot 3^4$, we can choose a wildly ramified Galois extension $E$ of $k$ with Galois group $\bZ/12\bZ$, and define a homomorphism $f:\Gal(k^s/k) \ra \Wabs$ by requiring that it factors through $\Gal(E/k)$ and sends a generator of $\Gal(E/k)$ to $w$. Analogous to above this yields an anisotropic maximal torus in the split inner form of $G$  that therefore transfers to $G$ and that does not split over any tamely ramified extension of $k$.

It remains to show that for $p=5$ there exists an inner form $G$ of every split group of type $E_6$ for which all tori split over a tamely ramified extension. Consider the inner form $G$ corresponding to the symbol $^1E^{16}_{6,2}$  in \cite{Tits-index}, which has the name $^3E_6$ in \cite{Tits}. This group has split rank two and the connected component of the centralizer $Z_G(S) \times_k k^s$ of a maximal split torus $S$ has type $A_2 \times A_2$, as can be read of from the index in Figure \ref{Figure-index-E} (\cite[page~58]{Tits-index}).

\begin{figure}[ht]
	\centering
	\begin{tikzpicture}
	
	\draw (0,0) -- (3.7,0);
	\draw (0,0) -- (4,0);
	\draw (2,0) -- (2,1);
	
	\draw[fill=black] (0,0) circle(.1);
	\draw[fill=black] (1,0) circle(.1);
	\draw[fill=black] (2,0) circle(.1);
	\draw (2,0) circle(0.2);
	\draw[fill=black] (2,1) circle(.1);
	\draw (2,1) circle(0.2);
	\draw[fill=black] (2,1) circle(.1);
	\draw[fill=black] (3,0) circle(.1);
	\draw[fill=black] (4,0) circle(.1);
	
%	\node at (-1,0) {$E_6$};	
	%	\node at (3,0) {$\cdots$};
	%	\node at (1.5, .3) {4};
	\node at (0,-.4) {$\alpha_1$};
	\node at (2.6,1) {$\alpha_2$};
	\node at (1,-.4) {$\alpha_3$};
	\node at (2,-.4) {$\alpha_4$};
	\node at (3,-.4) {$\alpha_5$};	
	\node at (4,-.4) {$\alpha_6$};
	
	\end{tikzpicture}
	\caption{Index for $^1E^{16}_{6,2}$ (\cite[page~58]{Tits-index})}	
	\label{Figure-index-E}	
\end{figure}

Suppose there exists a maximal torus $T$ of $G$ that does not split over a tamely ramified extension, and let $f \in \Hom(\Gal(k^s/k),\Wabs)$ denote a representative for the resulting cocycle when $T$ is viewed as a torus of the split inner form of $G$. Let $E$ be a Galois extension of $k$ such that $\Gal(k^s/E)$ is the kernel of $f \in \Hom(\Gal(k^s/k),\Wabs)$, and let $E^{ur}$ and $E^t$ denote the maximal unramified extension and maximal tamely ramified extension of $k$ contained in $E$, respectively. By the structure of local Galois groups, we have $\Gal(E/k) \simeq \Gal(E/E^t) \rtimes \left( \Gal(E^t/E^{ur}) \rtimes \Gal(E^{ur}) \right)$. Let $w_5$ be a generator of $f(\Gal(E/E^t))\simeq \bZ/5\bZ$.  Without loss of generality, $w_5=s_{\alpha_1}s_{\alpha_3}s_{\alpha_4}s_{\alpha_2}$ for the simple roots $\alpha_1, \alpha_2, \hdots, \alpha_6$ as provided by Figure  \ref{Figure-index-E}, i.e. $w_5$ corresponds to a Coxeter element in the subgroup $S_5 \triangleleft \Wabs$ generated by $\alpha_1, \alpha_2, \alpha_3$ and $\alpha_4$. 
Thus, if $w'$ is the image of a generator of $\Gal(E^t/E^{ur}) \subset \Gal(E/k)$ in $\Wabs$ (using the above choice of isomorphism $\Gal(E/k) \simeq \Gal(E/E^t) \rtimes \left( \Gal(E^t/E^{ur}) \rtimes \Gal(E^{ur}) \right)$ to obtain the embedding), it normalizes $\<w_5\>$, and we can write $w'= w'_1 \cdot w'_2$ with $w'_1 \in S_5 \triangleleft \Wabs$ and $w'_2 \in \Wabs$ such that $w'_2$ commutes with $w_5$. Note that the order of $w'_2w_5$ is a multiple of 5. By \cite[Table~9]{Carter} there exist precisely two conjugacy classes in $\Wabs$ whose elements have order divisible by 5, one is the conjugacy class of $w_5$, and the other is the conjugacy class of an element $w_{10}$ of order 10 that is the product of $w_5=s_{\alpha_1}s_{\alpha_3}s_{\alpha_4}s_{\alpha_2}$ and $w_2:=s_{\alpha_6}$.  %for the simple roots $\alpha_1, \alpha_2, \hdots, \alpha_6$ as provided by Figure \ref{Figure-index-E}. 
 Changing the decomposition $w'=w_1'\cdot w_2'$ if necessary (by multiplying $w_1'$ and $w_2'$ with suitable powers of $w_5$), we may assume that $w_2'$ is either trivial or has order two. If $w_2'$ has order two, then there exists $g\in \Wabs$ such that $gw_5w_2g^{-1}=w_5w_2'$, and raising the equation to the second and fifth power, we deduce that $g \in \Cent_W(w_5)$ and $w_2'=gw_2g^{-1}$. By \cite[Table~9]{Carter}, the conjugacy class of $w_5$ in $\Wabs$ contains $5184=2^6 \cdot 3^4$ elements, hence $\abs{\Cent_W(w_5)}=2 \cdot 5$, which implies that $\Cent_W(w_5)=\<w_5w_2\>$. Thus $w_2'=w_2$ in this case. Similarly, if $w''$ denotes the image of a generator of $\Gal(E^{ur}/k) \subset \Gal(E/k)$ in $\Wabs$, then $w''=w''_1\cdot w''_2$ with $w''_1 \in S_5$ and $w''_2 \in \{1, w_2\}$. 
 
Hence, we have
\begin{equation}
	 \< w_5 \> \subset  f(\Gal(k^s/k)) \subset S_5 \times \< w_2\> .
\end{equation}

Note that the subgroup $S_5$ of $\Wabs$ preservers the subspace $V_1$ of $X_*(T_{k^s}) \otimes_{\bZ} \bR$ spanned by $\check\alpha_1, \check\alpha_2, \check\alpha_3, \check\alpha_4$, and acts trivially on the subspace $V_2$ of $X_*(T_{k^s}) \otimes \bR$ spanned by $\check\omega_{\alpha_5}$ and the subspace $V_3$ of $X_*(T_{k^s}) \otimes \bR$ spanned by $\check\alpha_6$. In particular, the element $w_5$ acts on the subspace $V_1$ via a Coxeter element of $S_5$, and acts trivially on the subspace $V_2 \oplus V_3$. On the other hand, $w_2$ acts trivially on $V_1$ and $V_2$ and acts via the Coxeter element $-1$ of the subgroup $S_2 \triangleleft \Wabs$ generated by $s_{\alpha_6}$ on the subspace $V_3$. Therefore we have two possibilities for $\left(X_*(T_{k^s}) \otimes_{\bZ} \bR\right)^{\Gal(k^s/k)}=\left(V_1 \oplus V_2 \oplus V_3\right)^{\Gal(k^s/k)}$, which we treat separately. 

In the first case, $\left(X_*(T_{k^s}) \otimes_{\bZ} \bR\right)^{\Gal(k^s/k)}=\bR \cdot \check\omega_{\alpha_5} \oplus \bR \cdot \check\alpha_6=\bR \cdot \check\omega_{\alpha_5} \oplus \bR \cdot \check\omega_{\alpha_6}$, and hence $T$ has split rank two. If we denote by $S_T$ the maximal split subtorus of $T$, then by \cite[2.14~Lemma]{Steinberg-torsion} the reductive group $Z_G(S_T) \times_k k^s$ has root system of type $A_4$ (spanned by $\alpha_1, \alpha_2, \alpha_3, \alpha_4$), which contradicts that $S_T$ is conjugate to $S$ and $Z_G(S) \times_k k^s$ has root system of type $A_2 \times A_2$.

In the second case, $\left(X_*(T_{k^s}) \otimes_{\bZ} \bR\right)^{\Gal(k^s/k)}=\bR \cdot \check\omega_{\alpha_5}$, hence $T$ has split rank one. If we denote by $S_T$ the maximal split subtorus of $T$, then by \cite[2.14~Lemma]{Steinberg-torsion} the reductive group $Z_G(S_T) \times_k k^s$ has root system of type $A_4 \times A_1$ (spanned by $\alpha_1, \alpha_2, \alpha_3, \alpha_4$ and $\alpha_6$) and is supposed to contain a group isomorphic to $Z_G(S) \times_k k^s$ as a Levi subgroup. This is impossible since the latter has a root system of type $A_2 \times A_2$. 

Thus we obtain a contradiction in all cases and therefore every torus of $G$ splits over a tamely ramified extension. 
\qed

Using the result for absolutely simple groups, we can easily deduce the result for arbitrary reductive groups. 

\begin{Cor} \label{Cor-tame-tori}
	Let $G$ be a (connected) reductive group defined over a non-archimedean local field $k$ of residual characteristic $p$. 
	
	If $G$ splits over a tamely ramified field extension of $k$ and $p \nmid \abs{\Wabs}$, then every torus of $G$ splits over a tamely ramified field extension of $k$.

 If $\Phi(G)$ does not contain a component of type $A_n$ for some integer $n \geq 2$, $D_l$ for some prime $l \geq 4$ or $E_6$, then the reverse implication holds as well. Otherwise, the reverse direction might require slight modifications of the conditions on $p$ as described in Theorem \ref{Thm-tame-tori-less-easy} Part \ref{Prop-Part-2}.
  \end{Cor}

\Proof \\
Suppose $G$ splits over a tamely ramified Galois extension $E$ of $k$, and $p \nmid \abs{\Wabs}$. Let $T$ be a maximal torus of $G$, $T'$ a maximal split torus of $G_E$, and $W_{T'}=N_{G_E}(T')/T'$. Then $T_E$ is a maximal torus of $G_E$, and its isomorphism class is determined by a class of a cocycle in the image of $H^1(E,N_{G_E}(T')) \ra H^1(E,W_{T'}) \simeq \Hom(E,\Wabs)/\mathord\sim$. Since $p \nmid \Wabs$, the class becomes trivial over a tamely ramified Galois extension $E'$ of $E$, and hence $T$ splits over a tamely ramified extension of $k$.

In order to prove the reverse direction, suppose $p \mid \abs{\Wabs}$, because the case of $G$ not splitting over a tamely ramified extension is obvious. By taking the preimage, it suffices to exhibit a non-tame torus in the case that $G$ is adjoint. Then $G$ is a finite product of restrictions of scalars $\prod_{1 \leq i \leq r} \Res_{E_i/k} G_i$ for some finite Galois extensions $E_i$ of $k$ and $G_i$ absolutely simple reductive groups over $E_i$. If $p \mid \abs{\Wabs}$, the prime $p$ divides the order of the Weyl group of ${G_j} \times_{E_j} k^s$ for some $1 \leq j \leq r$.  If $\Phi(G_j)$ is not of type $A_n$ for some integer $n$ with $p \nmid n$, $D_p$ if $p \geq 4$, or $E_6$ if $p=5$, then by Theorem \ref{Thm-tame-tori-less-easy} there exists a torus $T_j$ in $G_j$ that does not split over tamely ramified extension of $k$. Thus we obtain a torus $\Res_{E_j/k} T_j$ of $G$ that does not split over a tamely ramified extension of $k$. \qed

\section{Existence of good elements} \label{section-good-elements}
Let $T$ be a maximal torus defined over $k$ and let $E$ be some finite separable extension of $k$. Following Moy--Prasad (\cite{MP1, MP2}) we define a filtration of $T(E)$ index by  $r \in \bR_{\geq 0}$ as follows. We let $T(E)_0$ denote the $\cO_E$-points of the identity component of the Néron lft-model of $T_E$, which is a subgroup of finite index in the maximal bounded subgroup of $T(E)$. For $r>0$ we set
$$ T(E)_r=\left\{t \in T(E)_0 \, | \, \val(\chi(t)-1) \geq r \, \, \forall \, \chi \in X^*(T_{k^s}) \right\} .$$

If $T$ splits over a tamely ramified extension of $E$, then for $r>0$ 
$$ T(E)_r=\left\{t \in T(E) \, | \, \val(\chi(t)-1) \geq r \, \, \forall \,  \chi \in X^*(T_{k^s}) \right\} .$$
[This can be seen as follows: If $T$ is split of rank $n$ over $E$, then the claim follows from the observation that the connected component of the Néron lft-model is isomorphic to $\bG_m^{n}$ (over $\cO_E$). If $T$ splits over a tamely ramified Galois extension $E'$ of 
$E$, then by \cite[4.6.2~Proposition]{Yu2}, we have $T(E)_r=(T(E')_r)^{\Gal(E'/E)}=T(E')_r \cap T(E)$, which yields the desired result using the split case.]

Similarly, we have a Moy--Prasad filtration of the Lie algebra $\ft(E):=\Lie(T)(E)$ given by
$$ \ft(E)_r=\left\{X \in \ft(E) \, | \, \val(d\chi(X)) \geq r \, \, \forall \,  \chi \in X^*(T_{k^s}) \right\} $$
for $r \in \bR$.

As commonly done, we set $T(E)_{r+}=\bigcup_{s>r}T(E)_s$ and $\ft(E)_{r+}=\bigcup_{s>r} \ft(E)_s$, and we sometimes abbreviate $\ft(k)_r$ by $\ft_r$ and write $\ft=\Lie(T)(k)$.

Note that different authors use different conventions for the indexing of the Moy--Prasad filtration. Here we have chosen our indices in such a way that if $E/k$ is tamely ramified Galois, then $\ft_r=\ft(E)_r \cap \ft$ (for $r \in \bR$) and $T(k)_r=T(E)_r \cap T(k)$ (for $r>0$).

Let us recall the definition of good elements in the Lie algebra introduced by \cite[Definition~2.2.4]{Adler}.
\begin{Def}
	A semisimple element $X \in \fg$ is called \textit{good} if there exists a maximal torus $T \subset G$ that splits over a finite tamely ramified Galois extension $E$ of $k$ and a real number $r$ such that $X \in \ft(E)_{r}-\ft(E)_{r+}$  and for every $\alpha \in \Phi(G,T_{E})$, we have $d\alpha(X)=0$ or $\val(d\alpha(X))=r$.
\end{Def}

These elements play an important role in the representation theory of $p$-adic groups and harmonic analysis on them, and authors commonly assume that all tori split over a tamely ramified extension and that every coset of $\ft_r/\ft_{r+}$ contains a good semisimple element. We have seen in Section \ref{section-tame-tori}  that all tori of an absolutely simple reductive group split over tamely ramified extension if $p \nmid \abs{\Wabs}$. In order to show that this also guarantees the existence of good semisimple elements, we will first draw some consequences from the condition that  $p \nmid \abs{\Wabs}$.

%Note that our definition is equivalent to  \cite[Definition~2.2.4]{Adler}, because the Moy--Prasad filtration behaves well with respect to tamely ramified base change, see \cite[Proposition~1.4.1]{Adler}.	

\begin{Lemma}
	\label{Lemma-implications-on-p}
	Let $G$ be a (connected) reductive group over $k$. Suppose $ p \nmid \abs{\Wabs}$. 
	Then 
	\begin{enumerate}[label=(\arabic*),ref=\arabic*]
		\item \label{item-bad-prime} $\bZ\Phi/\bZ\Phi_0$ is $p$-torsion free for every closed subsystem $\Phi_0$ of $\Phi:=\Phi(G)$, i.e. $p$ is not a bad prime for $\Phi(G)$ (in the sense of \cite[4.1]{Springer-Steinberg}),
		\item \label{item-fundamental-group} $X^*(T_{k^s})/\bZ[\alpha \, | \, \alpha \in \Phi]$ is $p$-torsion free,
		\item \label{item-index-of-connection} $p$ does not divide the index of connection of $\check \Phi_0$ for any sub-root system $\Phi_0$ of $\Phi$ that is generated by a proper subset of a basis of $\Phi$.
	\end{enumerate}
\end{Lemma}
\Proof\\
\eqref{item-bad-prime} It suffices to consider the case that $\Phi=\Phi(G)$ is irreducible. By \cite[4.3]{Springer-Steinberg}, the root system of type $A_n \, (n \geq 1)$ has no bad primes, $B_n \, (n\geq 2)$, $C_n \, (n\geq2)$ and $D_n \, (n\geq 4)$ have 2 as a bad prime, $E_6,E_7, F_4$ and $G_2$ have 2 and 3 as bad primes and $E_8$ has 2, 3 and 5 as bad primes. Comparing this with Table \ref{table-weyl-group}, we see that all bad primes divide $\abs{\Wabs}$.

\begin{table}[h]\footnotesize
	\begin{tabular}{|c|c|c|c|c|c|c|c|c|}
		\hline type  & 			$A_n \, (n \geq 1)$     & $B_n, C_n \, (n \geq 2)$  &  $D_n \, (n \geq 4)$              & $E_6$                 &  $E_7$                      &  $E_8$ & $F_4 $ & $G_2$  \\ 
		\hline bad primes &  -   & 2   &  2 & 2, 3 &  2, 3 & 2, 3, 5 &  2, 3 &  2,3 \\ 
		\hline index of connection &  $n+1$  & 2   &  4 & 3 &  2 & 1 &  1 & 1 \\ 
		\hline 
	\end{tabular} 
	\caption{Bad primes and index of connection of irreducible root systems}
	\label{table-index-of-connection}
\end{table}

\eqref{item-fundamental-group}
Recall that by \cite[VI.2~Proposition~7]{Bourbaki-4-6} the index of connection divides $\abs{\Wabs}$. (For the reader's convenience the index of connection is shown in Table \ref{table-index-of-connection}.)
%Hence $X_*(T_{k^s})/\bZ[H_\alpha \, | \, \alpha \in \Phi]$ is $p$-torsion free, because the order of the torsion subgroup of $X_*(T_{k^s})/\bZ[H_\alpha \, | \, \alpha \in \Phi]$ divides the index of connection, which is coprime to $p$.
Hence $X^*(T_{k^s})/\bZ[\alpha \, | \, \alpha \in \Phi]$ is $p$-torsion free, because the order of the torsion subgroup of $X^*(T_{k^s})/\bZ[\alpha \, | \, \alpha \in \Phi]$ divides the index of connection, which is coprime to $p$.

\eqref{item-index-of-connection} Note that the Weyl group of $\check\Phi$ and the Weyl group of $\Phi$ are isomorphic, and the Weyl group of $\check\Phi_0$ is a subgroup of $\check\Phi$, because a basis of $\check \Phi_0$ can be extended to a basis of $\check\Phi$ by \cite[VI.1, Proposition~4 and Proposition~24]{Bourbaki-4-6}. Hence $p$ does not divide the order of the Weyl group of $\check\Phi_0$ and therefore does not divide the index of connection of $\check\Phi_0$ either. 
\qed

The following theorem about the existence of good elements provides an alternative to correcting the proof of the existence of good elements in every Moy--Prasad filtration coset of the Lie algebra of a tamely ramified maximal torus in \cite[Lemma~5.12]{Adler-Roche} and, at the same time, uses much milder hypotheses than those of \cite{Adler-Roche} (for general reductive groups). For general reductive groups Adler and Roche (\cite{Adler-Roche}) assumed that $p>c(G)$, while we only assume $p \nmid \abs{\Wabs}$. For the precise definition of the constant $c(G)$ see \cite[p.~453]{Adler-Roche}. To give an example, for $E_6, E_7$ and $E_8$, the constant $c(G)$ is $\leq 113$, 373, and 1291, respectively (\cite[Table, p.~452]{Adler-Roche}), while the condition  $p \nmid \abs{\Wabs}$ is equivalent to $p$ being larger than 5, 7 and 7, respectively.

\begin{Thm} \label{Thm-good}
	Let $G$ be a reductive group defined over a non-archimedean local field $k$ of residual characteristic $p$. Suppose that $p \nmid \abs{\Wabs}$ and $G$ splits over a tamely ramified extension of $k$. Then for every maximal torus $T$ of $G$, and any $r \in \bR$, every coset of $\ft_r/\ft_{r+}$ contains a good element. %(Recall that $\ft=\Lie(T)(k)$.)
\end{Thm}

\Proof \\
Let $X + \ft_{r+}$ be a coset of  $\ft_r/\ft_{r+}$, i.e. pick $X \in \ft_r$. Let $\Phi_0 \subset \Phi:=\Phi(G,T_{k^s})$ be the collection of roots $\alpha$ for which $\val(d\alpha(X))>r$. If $\Phi_0 = \Phi$, then $X \subset \ft_{r+}$, because $X^*(T_{k^s})/\bZ[\alpha \, | \, \alpha \in \Phi]$ is $p$-torsion free by Lemma \ref{Lemma-implications-on-p}\eqref{item-fundamental-group}. In this case $X+\ft_{r+}=\ft_{r+}$ contains the good element 0. Thus, for the remainder of the proof we assume that $\Phi_0 \neq \Phi$.  Note that $\Phi_0$ is a closed subsystem of $\Phi$ (i.e. $\bZ \Phi_0 \cap \Phi = \Phi_0$). Since $p \nmid \abs{\Wabs}$, we have by Lemma \ref{Lemma-implications-on-p}\eqref{item-bad-prime} that $\bZ\Phi/\bZ\Phi_0$ is $p$-torsion free. Hence we obtain $\Phi_0 = \bQ\Phi_0 \cap \Phi$. Moreover, since $X$ and $T$ are defined over $k$, the set $\Phi_0$ is stable under the action of the absolute Galois group $\Gal(k^s/k)$. 
Let $\Delta_0$ be a basis for $\Phi_0$. For $\alpha \in \Delta_0$ we denote by ${\check\omega^0_{\alpha}}$ the element of $\bQ\check\Phi_0$ that satisfies $\<\check\omega^0_{\alpha},\beta\>=\delta_{\alpha,\beta}$ for all $\beta \in \Delta_0$.  Since  by Lemma \ref{Lemma-implications-on-p}\eqref{item-index-of-connection} the prime $p$ does not divide the index of connection of $\check \Phi_0$, i.e. $p \nmid \abs{\bZ[\check\omega^0_{\alpha}  \, | \, \alpha \in \Delta_0]/\bZ[\check\alpha \, | \, \alpha \in \Delta_0]}$, we deduce that $\check\omega^0_{\alpha} \in \bZ_{(p)}\Phi_0$. Hence, we can write $\check\omega^0_\alpha = \sum_{\beta \in \Delta_0} n_\beta \check\beta$ with $n_\beta \in\bZ_{(p)}$, and define $H_{\check\omega^0_{\alpha}}=\sum_{\beta \in \Delta_0} \ov n_\beta H_\beta \in \ft_{k^s}$, where $\ov n_\beta$ denotes the image of $n_\beta$ in $\cO_{k^s}$ (for $\beta \in \Delta_0$). Then $d\beta(H_{\check\omega^0_{\alpha}})=\delta_{\alpha, \beta}$ for $\beta \in \Delta_0$ and $d\beta(H_{\check\omega^0_{\alpha}}) \in \cO_{k^s}$ for all $\beta \in \Phi$. 

 Let $Y_1 \subset \ft_{k^s}=\Lie(T)(k^s)$ be the $k^s$-subspace spanned by $\{H_\alpha=d\check\alpha(1) \, | \, \alpha \in \Phi_0\}$ (or, equivalently, by $\{H_{\check\omega^0_{\alpha}} \, | \, \alpha \in \Phi_0\}$), and define $Y_2  \subset \ft_{k^s}$ by 
$$Y_2=\left\{ Z \in \Lie(T)({k^s}) \, | \, d\alpha(Z)=0 \, \, \forall \, \alpha \in \Phi_0 \right\}.$$
Then, $\ft_{k^s}=Y_1 \oplus Y_2$, and by the above observations about $\Phi_0$, the subspaces $Y_1$ and $Y_2$ are $\Gal(k^s/k)$-stable.
%, and if $H_\alpha \in Y$, then $\check\alpha \in \check\Phi_0$. 
Write $X=X_1+X_2$ with $X_1 \in Y_1$, $X_2 \in Y_2$. Since $Y_1$ and $Y_2$ are $\Gal(k^s/k)$-stable, $X_1$ and $X_2$ lie in $\ft$. We will show that $X_2$ is a good element in $X+\ft_{r+}$.

Note that for $\alpha \in \Phi_0$, we have
$d\alpha(X_2)=0$. Thus 
$$ X_1=\sum_{\alpha \in \Delta_0} d\alpha(X) H_{\omega^0_{\alpha}} , $$
and hence for all $\beta \in \Phi$
$$ \val(d\beta(X_1))=\val\left(\sum_{\alpha \in \Delta_0} d\alpha(X) d\beta(H_{\omega^0_{\alpha}})\right) >r . $$
%As $T$ splits over a tamely ramified field extension (by Corollary \ref{Cor-tame-tori}), the Moy--Prasad filtration (normalized with respect to the valuation that extends the normalized valuation on $k$) behaves well with respected to base change (in the sense of \cite[Proposition~1.4.1]{Adler}, but normalized so that no scaling of the depth is necessary). 
 Since $X^*(T_{k^s})/\bZ[\alpha \, | \, \alpha \in \Phi]$ is $p$-torsion free by Lemma \ref{Lemma-implications-on-p}\eqref{item-fundamental-group}, we obtain $X_1 \in \ft_{r+}$, and therefore $X_2 \in X+\ft_{r+}$. Moreover, since $d\alpha(X_2)=0$ for $\alpha \in \Phi_0$ and $\val(d\alpha(X_2))=r$ for $\alpha \in \Phi-\Phi_0$, the semisimple element $X_2$ is a good element. \qed

\begin{Rem}	\label{Rem-weaker-conditions-on-p}
	We did not make full use of the assumption that $p \nmid \abs{\Wabs}$, e.g. when $T$ is tame it suffices to assume that $p$ is not a bad prime for $\Phi(G)$, the quotient $X^*(T_{k^s})/\bZ[\alpha \, | \, \alpha \in \Phi]$ is $p$-torsion free, and that $p$ does not divide the index of connection of $\check \Phi_0$ for any $\Gal(k^s/k)$-stable subroot system $\Phi_0$ of $\Phi$ that is generated by a proper subset of a basis of $\Phi$. 
	
	If $G$ is an absolutely simple reductive group all of whose tori split over a tamely ramified extension, but $p \mid \abs{\Wabs}$, then $G$ has to be a non-split inner form of type $A_n$ and $p \nmid n+1$, or $D_p$ and $p \geq 4$, or $E_6$ and $p=5$. Hence $p$ is not a bad prime for $G$, see Table \ref{table-index-of-connection} (\cite[4.3]{Springer-Steinberg}), and $X^*(T_{k^s})/\bZ[\alpha \, | \, \alpha \in \Phi]$ is $p$-torsion free, because it divides $n+1$ for $A_n$, 4 for $D_p$ and 3 for $E_6$ (see Table \ref{table-index-of-connection}). Hence one might be able to deduce the existence of good elements by investigating the possibilities of $\Gal(k^s/k)$-stable subroot systems in $\Phi(G)$ given the restraint that all tori are tame, which implies that the $\Gal(k^s/k)$-action is sufficiently non-trivial, compare Remark \ref{Remark-inner-forms}.
	Some of these special cases have also already been dealt with in \cite[Proposition~5.4]{Adler-Roche}.
\end{Rem}

An analogous result can be obtained for good elements in the group $G(k)$, which we define as follows 
\begin{Def}
	Let $r \in \bR_{\geq0}$.
	A semisimple element $\gamma \in G(k)$ is called \textit{good of depth $r$} if there exists a maximal torus $T \subset G$ that splits over a finite tamely ramified Galois extension $E$ of $k$ such that $\gamma \in T(E)_{r}-T(E)_{r+}$  and for every $\alpha \in \Phi(G,T_{E})$, we have $\alpha(\gamma)=1$ or $\val(\alpha(\gamma)-1)=r$.
\end{Def}
Note that the definition of Adler and Spice (\cite[Definition~6.1]{Adler-Spice}) only assumes that the torus $T$ is ``tame-modulo-center''. However, our definition is more analogous to the Lie algebra case and has also been used in applications. %some applications, see e.g. \cite{He-Kim}.

\begin{Thm} \label{Thm-good-group}
	Let $G$ be a reductive group defined over a non-archimedean local field $k$ of residual characteristic $p$. Suppose that $p \nmid \abs{\Wabs}$ and $G$ splits over a tamely ramified extension of $k$. Then, for every maximal torus $T$ of $G$, and any $r \in \bR_{>0}$, every non-identity coset of $T(k)_r/T(k)_{r+}$ contains a good element of depth $r$. 
\end{Thm}
The proof of this theorem follows the spirit of the proof of Theorem \ref{Thm-good}. However, since some adjustments are necessary to pass from the Lie algebra to the group setting, we provide the details below.

\textbf{Proof of Theorem \ref{Thm-good-group}.}\\
Let $\gamma \in T(k)_r-T(k)_{r+}$. We consider the coset $\gamma \cdot T(k)_{r+}$. Let $\Phi_0 \subset \Phi:=\Phi(G,T_{k^s})$ be the collection of roots $\alpha$ for which $\val(\alpha(\gamma)-1)>r$. Then $\Phi_0$ is a closed subsystem of $\Phi$, and by Lemma \ref{Lemma-implications-on-p}\eqref{item-bad-prime} our assumption $p \nmid \abs{\Wabs}$ implies that $\bZ\Phi/\bZ\Phi_0$ is $p$-torsion free. Since $r>0$, we have $\chi(\gamma) \in 1 + \cP$ and $\val(\chi(\gamma)-1)=\val(\chi^n(\gamma)-1)$ for all characters $\chi \in X^*(T_{k^s})$ and positive integers $n$ coprime to $p$.
This implies that $\Phi_0 = \bQ\Phi_0 \cap \Phi$. Moreover, since $\gamma$ and $T$ are defined over $k$, the set $\Phi_0$ is $\Gal(k^s/k)$-stable. 
Let $\Delta_0$ be a basis for $\Phi_0$, and denote by $\check \omega^0_{\alpha}$ the element of $\bQ\check\Phi_0$ that satisfies $\<\check\omega^0_{\alpha},\beta\>=\delta_{\alpha,\beta}$ for all $\beta \in \Delta_0$.
Since  by Lemma \ref{Lemma-implications-on-p}\eqref{item-index-of-connection} the prime $p$ does not divide the index of connection of $\check \Phi_0$, for every $\alpha \in \Delta_0$ there exists an integer $n_\alpha$ coprime to $p$ such that $(\check \omega_\alpha^0)^{n_\alpha} \in \bZ\check\Phi_0 \subset X_*(T_{k^s})$. We define
$$ \gamma_1 = \prod_{\alpha \in \Delta_0} (\check \omega_\alpha^0)^{n_\alpha} ((\alpha(\gamma))^{1/n_\alpha}), $$
where $(\alpha(\gamma))^{1/n_\alpha}$ denotes the unique element $x$ in $1+\cP_{k^s}$ satisfying $x^{n_\alpha}=\alpha(\gamma)$.
We will show that $\gamma_2:=\gamma\gamma_1^{-1}$ is a good element of depth $r$ in $\gamma \cdot T(k)_{r+}$.

For $\beta \in \Phi$, we have 
\begin{eqnarray} \beta(\gamma_1)&=&\prod_{\alpha \in \Delta_0} \beta\left((\check \omega_\alpha^0)^{n_\alpha} ((\alpha(\gamma))^{1/n_\alpha})\right)
	=\prod_{\alpha \in \Delta_0}\left((\alpha(\gamma))^{1/n_\alpha}\right)^{n_\alpha\<\check \omega_\alpha^0,\beta\>} \notag 
	=\prod_{\alpha \in \Delta_0}\alpha(\gamma)^{\<\check \omega_\alpha^0,\beta\>}.
\end{eqnarray}
Note that since $p$ does not divide the index of connection of $\check \Phi_0$, we have $\<\check \omega_\alpha^0,\beta\> \in \bZ_{(p)}$, and $\alpha(\gamma)^{\<\check \omega_\alpha^0,\beta\>}$ makes sense inside $1+\cP_{k^s}$ (analogously defined as $\alpha(\gamma)^{1/n_\alpha}$ above).
Hence 
\begin{equation} \label{eqn-val-r} \val(\beta(\gamma_1)-1)=\val\left(\prod_{\alpha \in \Delta_0}\alpha(\gamma)^{\<\check \omega_\alpha^0,\beta\>}-1 \right)\geq\min\left\{\val(\alpha(\gamma)-1) \, | \, \alpha \in \Delta_0 \right\}>r, 
\end{equation}
and therefore $\gamma_1 \in T(k^s)_{r+}$, because $X^*(T_{k^s})/\bZ[\beta \, | \, \beta \in \Phi]$ is $p$-torsion free (Lemma \ref{Lemma-implications-on-p}\eqref{item-fundamental-group}). 
In addition, if $\beta \in \Phi_0$, then
\begin{equation} \label{eqn-val-1} \beta(\gamma_1)=\prod_{\alpha \in \Delta_0}\alpha(\gamma)^{\delta_{\alpha,\beta}}=\beta(\gamma) \quad \text{ and hence} \quad \beta(\gamma_2)=1 . \end{equation}
Therefore $\gamma_2$ is contained in the subgroup 
$$T_2:=\{t \in T(k^s) \, | \, \alpha(t)=1 \, \forall \, \alpha \in \Phi_0\}$$
 of $T(k^s)$, and $\gamma_1$ is by definition contained in the subgroup 
 $$T_1 :=\< \check\chi(1+\cP_{k^s})) \, | \, \check\chi \in \bQ\check\Phi_0\cap X_*(T_{k^s}) \>  \subset T(k^s) .$$
  Since $\Phi_0$ is preserved by the $\Gal(k^s/k)$-action, the subgroups $T_1$ and $T_2$ are $\Gal(k^s/k)$-stable. We claim that $T_1\cap T_2-\{1\}$ is empty. Suppose not, and let $t \in T_1\cap T_2-\{1\}$. Let $\check \chi_1, \hdots, \check \chi_n$ be a $\bZ$-basis for $\bQ\check\Phi_0\cap X_*(T_{k^s})$. Then $t= \check\chi_1(t_1) \cdot \hdots \cdot \check \chi_n(t_n)$ for some $t_1, \hdots, t_n \in 1+\cP_{k^s}$, one of which is non-trivial. Assume without loss of generality that $t_1 \neq 1$, and let $\chi_1 \in \bQ\Phi_0$ such that $\<\check\chi_i,\chi_1\>=\delta_{i,1}$ for $1 \leq i \leq n$. Then there exists an integer $N$ (coprime to $p$) such that $N{\chi_1}= \sum_{\alpha \in \Delta_0} a_\alpha \alpha$ with $a_\alpha \in \bZ$ for $\alpha \in \Delta_0$. 
  Hence $\chi:=N\chi_1 \in X^*(T_{k^s})$ with ${\chi}(t) = \prod_{\alpha \in \Delta_0} \alpha(t)^{a_\alpha}=1$ since $t \in T_2$. On the other hand, ${\chi}(t)=t_1^N \neq 1$, because $t_1 \in (1+\cP_{k^s}) -\{1\}$. This is a contradiction, and therefore $T_1 \cap T_2 = \{ 1\}$. Thus the factorization of $\gamma$ as $\gamma_1\gamma_2$ with $\gamma_1 \in T_1$ and $\gamma_2 \in T_2$ is unique, and hence $\gamma_1$ and $\gamma_2$ are contained in $T(k)$. 
  
  Since $T$ splits over a tamely ramified extension (by Corollary \ref{Cor-tame-tori}), we have $\gamma_1 \in T(k) \cap T(k^s)_{r+}=T(k)_{r+}$, hence $\gamma_2 \in \gamma \cdot T(k)_{r+}$. If $\alpha \in \Phi_0$, then $\alpha(\gamma_2)=1$ by Equation \eqref{eqn-val-1}, and
  if $\alpha \in \Phi -\Phi_0$, then we deduce from Inequality \eqref{eqn-val-r} that $\val(\alpha(\gamma_2)-1)=r$.
Thus $\gamma_2$ is a good element of depth $r$ in the coset $\gamma \cdot T(k)_{r+}$. \qed

\begin{Rem}
	As for Theorem \ref{Thm-good}, we did not make full use of the assumption that $p \nmid \abs{\Wabs}$. More precisely, Remark \ref{Rem-weaker-conditions-on-p} regarding weaker assumptions also applies to the proof of Theorem \ref{Thm-good-group}.
\end{Rem}

\bibliography{Jessicasbib}

\end{document}